\documentclass[12pt]{article}
\usepackage{amssymb,amsmath, amsthm, amscd,ifthen}

\usepackage{graphicx,url}

\usepackage[russian]{babel}   
\usepackage[utf8]{inputenc}  

\newcommand{\Natur}{{\mathbb N}}

\newtheorem{theorem}{Теорема}
\newtheorem{claim}{Утверждение}

\title{К статье о числе рёбер в индуцированных подграфах специального дистанционного графа}

\author{Ф.А. Пушняков \\ Московский Физико-Технический Институт \\
(Государственный Университет)
\\ email: filipp.pushnyakov@phystech.edu}

\begin{document} 

\maketitle

\section{Введение}

Рассмотрим последовательность графов $G_{n} = G_{n} (V_{n}, E_{n}) = G(n,3,1)$,  у которых $$V_{n} = \{x = (x_{1}, \dots, x_{n}) \; \vert \; x_{i} \in \{0, 1\},\; i = 1, \dots , n \; , \; x_{1} + \ldots + x_{n} = 3\},$$ $$E_{n} = \{ (x, y) \; \vert \; \langle x, y \rangle = 1\},$$
где через $\langle x, y \rangle$ обозначено скалярное произведение векторов $x$ и $y$. 
Иными словами, вершинами графа $G(n, 3, 1)$ являются $(0,1)$-векторы, скалярный квадрат которых равен трем. И эти вершины соединены ребром тогда и только тогда, когда скалярное произведение соответствующих векторов равно единице. Данное определение можно переформулировать в комбинаторных терминах. А именно, рассмотрим граф, вершинами которого являются всевозможные трехэлементные подмножества множества $\mathcal{R}_{n} = \{1, \dots, n\}$, причем ребро между такими вершинами проводится тогда и только тогда, когда соответствующие трехэлементные подмножества имеют ровно один общий элемент. Изучение данного графа обусловлено многими задачами комбинаторной геометрии, экстремальной комбинаторики, теории кодирования: например, задачей Нелсона--Эрдёша--Хадвигера о раскраске метрического пространства (см. \cite{Rai3}--\cite{Kiselev1}), проблемой Борсука о разбиении пространства на части меньшего диаметра (см. \cite{Rai3}--\cite{Rai5}, \cite{Bolt}--\cite{Rai2}), задачами о числах Рамсея (см. \cite{Ram}--\cite{Rai12}), задачами о кодах с одним запрещенным расстоянием (см. \cite{MS}--\cite{Bobu3}).
\par Напомним несколько свойств данного графа. Граф $G(n,3,1)$ является регулярным со степенью вершины $d_{n} = 3 \cdot C_{n-3}^{2}$. Очевидно, что $\vert V_{n} \vert = C_{n}^{3} \sim \frac{n^{3}}{6}$ при  $n \rightarrow \infty$. В силу регулярности рассматриваемого графа имеем $\vert E_{n} \vert = \frac{d_{n} \cdot \vert V_{n} \vert}{2} = \frac{3}{2} \cdot C_{n-3}^{2} \cdot C_{n}^{3} \sim \frac{n^{5}}{8}$ при $n \rightarrow \infty$. \par
Напомним, что \textit{независимым множеством} графа называется такое подмножество его вершин, что никакие две вершины подмножества не соединены ребром. \textit{Числом независимости} $\alpha(G)$ называется наибольшая мощность независимого множества.  Положим $\alpha_{n} = \alpha(G(n, 3, 1))$. Результат теоремы Ж. Надя (см. \cite{Nagy}) отвечает на вопрос о числе независимости графа $G(n, 3, 1)$. А именно, $\alpha_{n} \sim n$ при $n \rightarrow \infty$. Более того, из доказательства теоремы Ж. Надя можно сделать вывод о структуре независимого множества в рассматриваемом графе. Для описания этой структуры введем дополнительные обозначения.  Пусть $W \subseteq V_{n}$. Будем говорить, что $W$ является множеством вершин \textit{первого типа}, если $\vert W \vert \geq 3$ и существуют такие $i, j \in \mathcal{R}_{n}$, что для любой вершины $w \in W$ выполнено $i, j \in w$; далее, $W$ является множеством вершин \textit{второго типа}, если $\vert W \vert \geq 2$ и существуют такие $i, \; j, \; k, \; t \in \mathcal{R}_{n}$, что для любой вершины $w \in W$ выполнено $w \subset \{i, j, k, t\} $; наконец, $W$ является множеством вершин \textit{третьего типа}, если для любых $w_{1}, w_{2} \in W$ выполнено соотношение $w_{1} \cap w_{2} = \emptyset$. Более того, \textit{носителем} множества вершин назовем объединение всех вершин данного множества. Тогда имеет место следующее утверждение.
\begin{claim} Любое независимое множество $U \subseteq V_{n}$ можно представить в виде объединения $$U = \left( \cup_{i \in \mathcal{I}} A_{i} \right) \cup \left( \cup_{j \in \mathcal{J}} B_{j} \right) \cup \left( \cup_{k \in \mathcal{K}} C_{k} \right),$$ где $A_{i}$ -- множество вершин первого типа, $B_{j}$ -- множество вершин второго типа, $C_{k}$ -- множество вершин третьего типа, $i \in \mathcal{I}, \; j \in \mathcal{J},\;  k \in \mathcal{K}$, и носители всех упомянутых множеств попарно не пересекаются.
\end{claim}
Мы не доказываем данное утверждение, так как оно мгновенно следует из доказательства теоремы Ж. Надя (см. \cite{Nagy}).
\par
Обозначим через $r(W)$ количество рёбер графа $G$ на множестве $W \subseteq V_{n}$. Иными словами, $$r(W) = \vert \{(x, y) \in E(G) \; \vert \; x \in W, \; y \in W\} \vert \; .$$ Также положим  $$r(l(n)) = \min_{\vert W \vert = l(n), \; W \subseteq V_{n}} r(W) \; .$$
Заметим, что если $l(n) \leq \alpha_{n}$, то $r(l(n)) = 0$ и обсуждать нечего. Если же $l(n) > \alpha_{n}$, то, очевидно, в любом $W \subseteq V_{n}$ мощности $l(n)$ непременно найдутся рёбра. Возникает интересный вопрос об изучении величины $r(l(n))$. Оценки, полученные в работах \cite{Rai11}--\cite{Rai12}, достаточно слабые, поэтому появились работы \cite{Pushnyakov}--\cite{Pushnyakov2}, в которых приведено практически полное исследование величины $r(l(n))$. А именно, в работе \cite{Pushnyakov} была доказана следующая теорема (читая формулировку, важно помнить, что $ n \sim \alpha_n $).

\begin{theorem} \label{t1}
Имеют место четыре случая:
\begin{enumerate}
\item Пусть функции $f: \Natur \rightarrow \Natur, \; g: \Natur \rightarrow \Natur$ таковы, что выполнено $n = o(f)$ и $g = o(n^{2})$ при $n \rightarrow \infty$. Пусть функция $l: \Natur \rightarrow \Natur$ такова, что для любого $n \in \Natur$  выполнена цепочка неравенств $f(n) \leq l(n) \leq g(n)$. Тогда $r(l(n)) \sim \frac{l(n)^{2}}{2 \alpha_{n}}$ при $n \rightarrow \infty$.
\item Пусть функция $l: \Natur \rightarrow \Natur$  такова, что существуют константы $C_{1}, \; C_{2}$, с которыми для каждого $n \in \Natur$ выполнена цепочка неравенств $C_{1} \cdot n^{2} \leq f(n) \leq C_{2} \cdot n^{2}$. Тогда $r(l(n)) \sim \frac{l(n)^{2}}{2 \alpha_{n}}$ при $n \rightarrow \infty$.
\item Пусть функции $f: \Natur \rightarrow \Natur, \; g: \Natur \rightarrow \Natur$ таковы, что выполнено $n^{2} = o(f)$ и $g = o(n^{3})$ при $n \rightarrow \infty$.  Пусть функция $l: \Natur \rightarrow \Natur$ такова, что для каждого $n \in \Natur$ выполнено $f(n) \leq l(n) \leq g(n)$. Тогда существуют такие функции $h_{1}: \Natur \rightarrow \Natur, \ h_{2}: \Natur \rightarrow \Natur$, что $h_{1} \sim \frac{l(n)^{2}}{\alpha_{n}}, \ h_{2} \sim \frac{5 l(n)^{2}}{\alpha_{n}}$  при $n \rightarrow \infty$ и для каждого $n \in \Natur$ выполнена цепочка неравенств $h_{1}(n) \leq r(l(n)) \leq h_{2}(n)$. При этом для выполнения нижней оценки требование $ g = o(n^3) $ не нужно.
\item Пусть функция $l: \Natur \rightarrow \Natur$  такова, что существует константа $C$, с которой выполнена цепочка неравенств $C \cdot n^{3} \leq l(n) \leq  C_{n}^{3}$. Пусть $c_{n} = 1-\frac{l(n)}{C_{n}^{3}}$. Тогда существует функция $f: \Natur \rightarrow \Natur$, такая, что $f(n) \sim n^{5} \left(\frac{1}{8} - \frac{c_{n}}{4} + \frac{c_{n}^{2}}{72}\right)$ при $n \rightarrow \infty$ и для каждого $n \in \Natur$ выполнено $r(l(n)) \geq f(n)$.
\end{enumerate}
\end{theorem}

Как можно заметить, оценки, полученные в пунктах 3-4 данной теоремы, не являются точными. Отметим в то же время, что нижняя оценка из пункта 3 верна и в условиях пункта 4. В работе \cite{Pushnyakov2} была улучшена оценка из пункта 3. А именно, была доказана следующая теорема.

\begin{theorem} \label{t3}
Пусть функция $l: \Natur \rightarrow \Natur$ такова, что $n^{2} = o(l)$ при $n \rightarrow \infty$. Тогда существует такая функция $h: \Natur \rightarrow \Natur$, что $h \sim \frac{3 l^{2}}{2n}$ при $n \rightarrow \infty$ и $r(l(n)) \geq h(n)$ для достаточно большого $n$. 
\end{theorem}

В данной работе автор устранил одну неточность в пункте 4 теоремы 1 и улучшил многие оценки из теоремы 1, а именно доказал следующую теорему.

\begin{theorem} \label{t2}
Имеют место четыре случая:
\begin{enumerate}
\item
Пусть дана произвольная функция $l: \Natur \rightarrow \Natur$ с ограничением $n = o(l)$. Тогда существует такая функция $h: \Natur \rightarrow \Natur$, что $h(n) \sim \frac{9 l(n)^{2}}{2 \alpha_{n}} $  при $n \rightarrow \infty$ и для каждого $n \in \Natur$ выполнено неравенство $r(l(n)) \leq h(n)$.
\item Пусть функция $l: \Natur \rightarrow \Natur$  такова, что существуют константа $C$ и функция $g: \Natur \rightarrow \Natur$ такая, что $n^{2} = o(g)$ при $n \rightarrow \infty$ и для каждого $n \in \Natur$ выполнена цепочка неравенств $C \cdot n^{3} \leq l(n) \leq  C_{n}^{3} - g(n)$. Пусть $c_{n} = 1-\frac{l(n)}{C_{n}^{3}}.$ 
Тогда существует такая функция $h : \Natur \rightarrow \Natur$, что $h = o\left(1\right)$ при $n \rightarrow \infty$ и для любого достаточно большого $n \in \Natur$ выполнено неравенство $r(l(n)) \geq \frac{n^{5}}{8} \left(1 - 2c_{n} + \frac{c_{n}^{2}}{3} \left(1 + h\left(n\right)\right) - \frac{10}{n} + \frac{20 \cdot c_{n}}{n} - \frac{10 \cdot c_{n}^{2}}{3 n} \left(1 + h\left(n\right)\right)\right).$
\item Пусть функция $l: \Natur \rightarrow \Natur$  такова, что существуют константы $B, \ C$ и функция $g: \Natur \rightarrow \Natur$ такая, что $n = o\left(g\right)$ при $n \rightarrow \infty$ и для каждого $n \in \Natur$ выполнено $g(n) \leq Bn^{2}$ и $C \cdot n^{3} \leq l(n) \leq  C_{n}^{3} - g(n)$. Пусть $c_{n} = 1-\frac{l(n)}{C_{n}^{3}}.$ 
Тогда существует такая функция $h : \Natur \rightarrow \Natur,$ что $h = o\left(1\right)$ при $n \rightarrow \infty$ и для любого достаточно большого $n \in \Natur$ выполнено неравенство $r(l(n)) \geq \frac{n^{5}}{8} \left(1 - 2 c_{n} + \frac{2 c_{n}^{2}}{9} \left(1 + h\left(n\right)\right) - \frac{10}{n} + \frac{20 \cdot c_{n}}{n} - \frac{20 \cdot c_{n}^{2}}{9 n} \left(1 + h\left(n\right)\right)\right).$
\item Пусть функция $l: \Natur \rightarrow \Natur$  такова, что существует такая константа $C$, что для каждого $n \in \Natur$ выполнена цепочка неравенств $C_{n}^{3} - Cn \leq l(n) \leq  C_{n}^{3}$. Пусть $c_{n} = 1-\frac{l(n)}{C_{n}^{3}}$. 
Тогда для любого $n \in \Natur$ выполнено неравенство $r(l(n)) \geq \frac{n^{5}}{8} \left(1 - 2 c_{n} - \frac{10}{n} + \frac{20 \cdot c_{n}}{n}\right).$

\end{enumerate}
\end{theorem}
Результат первого пункта теоремы \ref{t2} улучшает оценку, полученную в пункте 3 теоремы \ref{t1}. Тем не менее, эта оценка по-прежнему не является точной: величина $r(l(n))$ удовлетворяет следующей цепочке неравенств: 
$$\frac{3 l(n)^{2}}{2 \alpha_{n}} (1 + o(1)) \leq r(l(n)) \leq \frac{9 l(n)^{2}}{2 \alpha_{n}} (1 + o(1))$$
при $n \rightarrow \infty$. Между нижней и верхней оценками имеется зазор в 3 раза.
\par 
Результат второго и третьего пунктов теоремы \ref{t2} немного улучшает оценку, полученную в пункте 4 теоремы \ref{t1}. А пункт 4 как будто слабее аналогичного пункта из первой теоремы. Дело в том, что в доказательстве теоремы 1 была допущена неточность и здесь мы её устраняем.  К тому же, оценки из пунктов 2--4 теоремы 3 выглядят несколько иначе, чем аналогичная оценка из теоремы 1. Дело в том, что в формулировке и доказательстве пункта 4 теоремы 1 была допущена ошибка при переходе к асимптотикам, делающая результат некорректным при $c_{n} \rightarrow 0$. Поэтому вместо того, чтобы искать асимптотику функции $f$, автор нашёл её точное значение. Таким образом, теперь переход к пределу при $c_{n} \rightarrow 0$ абсолютно корректен.
\par Посмотрим ещё с несколько иной точки зрения на полученные результаты. Можно записать лучшую известную нам верхнюю оценку в виде 
\begin{equation}\label{formula1}
	r(l(n)) \le \frac{9l(n)^{2}}{2n} (1+o(1)) = \frac{n^{5}}{8} (1-2c_{n}+c_{n}^{2}) (1+o(1)),
\end{equation}
где $ c_n $ из формулировки теоремы 3. В таком же виде можно записать и нижнюю оценку из теоремы 2:
\begin{equation}\label{formula2}
	r(l(n)) \ge \frac{n^5}{8} \left(\frac{1}{3}-\frac{2}{3}c_{n}+\frac{1}{3} c_{n}^2\right) (1+o(1)).
\end{equation}

Как мы помним, в условиях пункта 4 теоремы 1 она верна, то есть верна она и в условиях пунктов 2--4 теоремы 3. Конечно, если $ c_{n} \rightarrow 0 $, то оценки пунктов 2--4 новой теоремы асимптотически совпадают с оценкой (\ref{formula1}) и в этом случае оценка (\ref{formula2}) им не конкурент. Однако в условиях теоремы 3 возможно и что $ c_{n} $ не стремится к нулю (хотя и не превосходит константы, строго меньшей единицы). В этом случае оценки из пунктов 2--4 становятся лучше, чем оценка из теоремы 2 при выполнении неравенства  
$$\frac{n^{5}}{8} \left(1 - 2 c_{n} - \frac{10}{n} + \frac{20 \cdot c_{n}}{n}\right) \geq \frac{3 \left( c_{n} C_{n}^{3}\right)^{2}}{2},$$
которое выполнено при $c_{n} \le 0.486 \dots$ и достаточно больших значениях $n \in \Natur$.

\par

\section{Доказательство теоремы 3} \label{sec:firstpage}
\subsection{Доказательство пункта 1}
Для доказательства верхней оценки необходимо для каждой функции $l(n)$, удовлетворяющей условию пункта 1 теоремы 3, и для каждого $n$ построить пример множества $W_{n}$ мощности $l(n)$, для которого величина $r(W_{n})$ оценивается сверху нужным образом. По-прежнему можно считать, что $n$ достаточно велико. \par 
Зафиксируем произвольную функцию $l$, удовлетворяющую условию пункта 1 теоремы 3, и число $n$. Возьмём наименьшее натуральное число $t_{n}$, с которым $C_{t_{n}}^{3} \cdot \left[ \frac{n}{t_{n}} \right] \ge l(n)$. Ясно, что $t_{n} \rightarrow \infty$ при $n \rightarrow \infty$. Значит, $t_{n} \sim \sqrt{\frac{6l}{n}}$ при $n \rightarrow \infty$. Положим
$$S_{1} = \left\{ 1, 2, \dots, t_{n} \right\},$$
$$S_{2} = \left\{t_{n} + 1, t_{n} + 2, \dots, 2 \cdot t_{n} \right\},$$
$$\dots$$
$$S_{\left[ \frac{n}{t_{n}} \right]} = \left\{ \left[ \frac{n}{t_{n}} \right] \cdot t_{n} - t_{n} + 1, \left[ \frac{n}{t_{n}} \right] \cdot t_{n} - t_{n} + 2, \dots, \left[ \frac{n}{t_{n}} \right] \cdot t_{n} \right\}.$$
Подчеркнём, что мощность каждого из множеств $S_{i},\ i=1, \dots, \left[ \frac{n}{t_{n}} \right],$ равна $t_{n}$. 
\par 
Для каждого $i=1, \dots, \left[ \frac{n}{t_{n}} \right]$ положим 
$$U_{i} = \bigcup_{x \in S_{i}, \ y \in S_{i} \ z \in S_{i}, \; x \neq y \neq z} \left\{ \left\{x, y, z\right\} \right\}.$$
Иными словами, множество $U_{i}$ --- это подмножество множества вершин графа $G(n, 3, 1),$ носители которых лежат во множестве $S_{i}$. \par 
Заметим, что в силу выбора величины $t_{n}$ данный выбор подмножеств корректен. 
То есть, для каждого множества $S_{i}$ существуют вершины графа $G(n, 3, 1)$ с носителем, лежащим в $S_{i}$.
Коль скоро $U_{i}$ --- это подмножество множества вершин графа $G(n, 3, 1)$, мы можем положить $$T_{i} = \left\{(v, w) \ \vert \ v \in U_{i}, \ w \in U_{i}, \ (v, w) \in E_{n}\right\}.$$

Иными словами, $T_{i}$ -- это множество рёбер графа $G(n, 3, 1)$, вершинами которых являются вершины из множества $U_{i}.$ Посчитаем мощности множеств $U_{i}$ и $T_{i}$. \par 
Легко видеть, что 
$$\vert U_{i} \vert = C_{\vert S_{i} \vert}^{3} = \frac{\vert S_{i} ^{3}\vert}{6} (1 + o(1)) = \frac{t_{n}^{3}}{6} (1 + o(1)).$$
Для оценки величины $\vert T_{i} \vert$ заметим, что подграф графа $G(n, 3, 1)$, индуцированный подмножеством вершин $U_{i}$, является регулярным. Положим $d_{i}$ -- степень вершины в подграфе графа $G(n, 3, 1)$, индуцированном подмножеством вершин $U_{i}$. Легко видеть, что 
$$d_{i} = 3 C_{\vert S_{i} \vert - 3}^{2} = \frac{3 \vert S_{i} \vert ^{2}}{2} (1 + o(1)) = \frac{3 t_{n}^{2} }{2} (1 + o(1)).$$
Стало быть, 
$$\vert T_{i} \vert = \frac{\vert U_{i}\vert \cdot d_{i}}{2} = \frac{t_{n}^{5} }{8} (1 + o(1)).$$
Перейдём к построению множества $W_{n}$, обладающего описанными в начале доказательства свойствами. Положим 
$$W_{n} = \bigcup_{i=1}^{\left[ \frac{n}{t_{n}} \right]} U_{i}.$$
Заметим, что для любых $i, \ j \leq \left[ \frac{n}{t_{n}} \right]$, множества вершин $U_{i}$ и $U_{j}$ не пересекаются. Таким образом, объединение выше на самом деле является дизъюнктным. 
\par Положим 
$$E(W_{n}) = \left\{ (u, v) \ \vert \ u \in W_{n}, \ v \in W_{n}, \ (u, w) \in E_{n} \right\}.$$
Иными словами, $E(W_{n})$ --- это множество рёбер в подграфе графа $G(n, 3, 1)$, индуцированным вершинами $W_{n}$. Поскольку множество $W_{n}$ является дизъюнктным объединением множеств вершин $U_{i}$, которые попарно не пересекаются, то и множество рёбер $E(W_{n})$ является дизъюнктным объединением множеств рёбер $T_{i}$: 
$$E(W_{n}) = \bigsqcup_{i=1}^{\left[ \frac{n}{t_{n}} \right]} T_{i}.$$
Оценим мощности множеств $W_{n}$ и $E(W_{n})$. Из вышесказанного следует, что 
$$\vert W_{n} \vert = \sum_{i=1}^{ \left[ \frac{n}{t_{n}} \right]} \vert U_{i} \vert = \left[ \frac{n}{t_{n}} \right] \cdot C_{t_{n}}^{3} \ge l(n).$$
Также ясно, что  $\vert W_{n} \vert = \frac{n \cdot t_{n}^{2}}{6}(1 + o(1))$ при $n \rightarrow \infty$. Далее,
$$\vert E\left(W_{n}\right) \vert = \sum_{i=1}^{\left[ \frac{n}{t_{n}} \right]} \vert T_{i} \vert = \left[ \frac{n}{t_{n}} \right] \cdot \vert T_{1} \vert \cdot \left(1 + o\left(1\right)\right) = \frac{n \cdot t_{n}^{4}}{8}\left(1 + o\left(1\right)\right).$$
Таким образом, 
$$\vert E\left(W_{n}\right) \vert = \frac{n \cdot t_{n}^{4}}{8}(1 + o\left(1\right)) = \frac{9}{2} \frac{\vert W_{n} \vert^{2}}{\alpha_{n}} (1+o\left(1\right)).$$
Стало быть, для каждой функции $l\left(n\right)$, удовлетворяющей условию пункта 1 теоремы, и для каждого $n$ мы построили пример множества $W_{n}$ мощности $l\left(n\right)$, для которого величина $r\left(W_{n}\right)$ оценивается сверху нужным образом. 

\subsection{Доказательство пунктов 2--4}
Для произвольной функции $l$, удовлетворяющей неравенствам $n^{2} = o\left(l\right)$ и $l \leq C_{n}^{3}$, и произвольного натурального числа $n$ положим $c_{n} = 1 - \frac{l(n)}{C_{n}^{3}}$.
Рассмотрим произвольное подмножество вершин $W \subseteq V_{n}$ мощности $l\left(n\right)$, положим $W_{1} = V_{n} \setminus W$. Ясно, что $\vert W_{1} \vert = c_{n} C_{n}^{3}$. Обозначим через $E\left(W_{1}\right)$ множество рёбер, концами которых являются вершины из $W_{1}$. Формально, 
$$E(W_{1}) = \left\{\left(x, y\right) \in E_{n} \; \vert \; x, y \in W_{1}\right\}.$$ 
Обозначим через $E_{1}$ множество рёбер, один конец которых принадлежит множеству $W$, а другой --- множеству $W_{1}$: 
$$E_{1} = \left\{\left(x, y\right) \in E_{n} \; \vert \; x \in W, \; y \in W_{1} \right\}.$$ 
С учетом введенных обозначений мы имеем $$E\left(W\right) = E_{n} \setminus \left( E\left(W_{1}\right) \sqcup E_{1} \right).$$ Тогда ясно, что
$$\vert E(W) \vert = \vert E_{n} \vert - \vert E(W_{1}) \vert - \vert E_{1} \vert.$$
Оценим сверху величину $\vert E(W_{1}) \vert + \vert E_{1} \vert$. В силу регулярности графа $G(n,3,1)$ имеем $$\vert E(W_{1}) \vert + \vert E_{1} \vert \leq d_{n} \cdot \vert W_{1}\vert.$$ Действительно, каждое ребро из множеств  $E(W_{1}) \cup E_{1}$ имеет одним из своих концов вершину из  $W_{1}$. Поэтому этих рёбер не больше, чем общее число рёбер, содержащих вершины из $W_{1}$. Данную оценку можно слегка уточнить. 
 
Заметим, что при таком подсчете дважды были посчитаны рёбра из $E(W_{1})$. Оценим мощность данного множества. Возможны три случая.
\begin{enumerate}
	\item \textit{Функция $l$ удовлетворяет условиям пункта 2 теоремы 3.}
	
	В данном случае мощность множества рёбер $E(W_{1})$ можно оценить с помощью результата теоремы 2. А именно, верна оценка:
	$$\vert E\left(W_{1}\right) \vert \geq \frac{{3} \vert W_{1} \vert^{2}}{2\alpha_{n}} \left(1+o\left(1\right)\right).$$
	Формально говоря, это означает, что существует такая функция $h : \Natur \rightarrow \Natur$, что $h = o\left(1\right)$ при $n \rightarrow \infty$ и для любого натурального $n$ выполнено неравенство 
	$$\vert E(W_{1}) \vert \geq \frac{{3} \vert W_{1} \vert^{2}}{2\alpha_{n}} \left(1+h\left(n\right)\right).$$
	Тогда $$\vert E(W_{1}) \vert + \vert E_{1} \vert \leq d_{n} \cdot \vert W_{1}\vert - \frac{3 \vert W_{1} \vert^{2}}{2\alpha_{n}}\left(1+h\left(n\right)\right).$$
	Положим 
\begin{equation*}
f\left(n\right) = \begin{cases}
0 &\text{\text{при } \; $n \leq 10$}\\
\frac{n^{5}}{8} \left(1 - 2c_{n} + \frac{c_{n}^{2}}{3} \left(1 + h\left(n\right)\right) - \frac{10}{n} + \frac{20 \cdot c_{n}}{n} - \frac{10 \cdot c_{n}^{2}}{3 n} \left(1 + h\left(n\right)\right)\right) &\text{\text{при } $n > 10$}
\end{cases}
\end{equation*}
Отметим, что константа 10 в определении функции $f$ не несёт значительной смысловой нагрузки --- она нужна лишь для упрощения некоторых выкладок ниже. \par
Заметим, что $f\left(n\right) \sim n^{5} \left(\frac{1}{8} - \frac{c_{n}}{4} + \frac{c_{n}^{2}}{24} \right)$ при $n \rightarrow \infty$.  Стало быть, осталось доказать, что $r\left(l\left(n\right)\right) \geq f\left(n\right)$ для любого достаточно большого натурального $n$.
	Для этого запишем следующую цепочку неравенств:
	$$r(l(n)) \geq \vert E(W) \vert \geq \frac{3}{2} C_{n-3}^{2} C_{n}^{3} - d_{n} \cdot \vert W_{1} \vert + \frac{3\vert W_{1} \vert^{2}}{2\alpha_{n}}\left(1+h\left(n\right)\right) \geq$$
$$\geq \frac{3}{2} C_{n-3}^{2} C_{n}^{3} - 3 \cdot C_{n-3}^{2} \cdot \vert W_{1} \vert + \frac{3\vert W_{1} \vert^{2}}{2n} \left(1 + h\left(n\right)\right) =$$
$$ =  \frac{3}{2} C_{n-3}^{2} C_{n}^{3} - 3 \cdot C_{n-3}^{2} \cdot c_{n} \cdot C_{n}^{3} + \frac{3 \cdot c_{n}^{2} \cdot \left( C_{n}^{3} \right) ^{2}}{2n} (1 + h(n)) = $$
$$ = \frac{3}{2}C_{n-3}^{2} C_{n}^{3} \left(1 - 2c_{n} + \frac{c_{n}^{2} \cdot C_{n}^{3}}{n \cdot C_{n-3}^{2}} \left(1 + h(n) \right)\right) = $$
$$ = \frac{n^{5}}{8} \left(1 - \frac{10}{n} + \frac{35}{n^{2}} - \frac{50}{n^{3}} + \frac{24}{n^{4}}\right) \cdot \left(1 - 2c_{n} + \frac{c_{n}^{2} \cdot \left( n - 1\right) \left( n - 2\right)}{3 \left(n - 3 \right) \left(n - 4 \right)} \left(1 + h \left( n \right) \right) \right) \geq $$
$$ \geq \frac{n^{5}}{8} \left(1 - \frac{10}{n} \right) \cdot \left(1 - 2c_{n} + \frac{c_{n}^{2}}{3} \left(1 + h(n) \right) \right) \geq f(n).$$


	\item \textit{Функция $l$ удовлетворяет условиям пункта 3 теоремы 3.}
	
	В данном случае воспользуемся результатами пунктов 1 и 2 теоремы 1. В силу этих результатов имеет место оценка: 
	$$\vert E\left(W_{1}\right) \vert \geq \frac{\vert W_{1} \vert^{2}}{\alpha_{n}} \left(1+o\left(1 \right)\right).$$
	Формально говоря, это означает, что существует такая функция $h : \Natur \rightarrow \Natur$, что $h = o\left(1\right)$ при $n \rightarrow \infty$ и для любого натурального $n$ выполнено неравенство 
	$$\vert E\left(W_{1}\right) \vert \geq \frac{\vert W_{1} \vert^{2}}{\alpha_{n}}\left(1+h\left(n\right)\right).$$
	Тогда 
	$$\vert E(W_{1}) \vert + \vert E_{1} \vert \leq d_{n} \cdot \vert W_{1}\vert - \frac{\vert W_{1} \vert^{2}}{\alpha_{n}}\left(1+h\left(n\right)\right).$$
	
	Положим 
\begin{equation*}
f\left(n\right) = \begin{cases}
0 &\text{\text{при } \; $n \leq 10$}\\
\frac{n^{5}}{8} \left(1 - 2 c_{n} + \frac{2 c_{n}^{2}}{9} \left(1 + h\left(n\right)\right) - \frac{10}{n} + \frac{20 \cdot c_{n}}{n} - \frac{20 \cdot c_{n}^{2}}{9 n} \left(1 + h\left(n\right)\right)\right) &\text{\text{при } $n > 10$}
\end{cases}
\end{equation*}

Заметим, что $f\left(n\right) \sim n^{5} \left(\frac{1}{8} - \frac{c_{n}}{4} + \frac{c_{n}^{2}}{36} \right)$ при $n \rightarrow \infty$. Стало быть, как и в предыдущем пункте, осталось доказать, что $r\left(l\left(n\right)\right) \geq f\left(n\right)$ для любого достаточно большого натурального $n$.
	
	Для этого запишем следующую цепочку неравенств:
	$$r(l(n)) \geq \vert E(W) \vert \geq \frac{3}{2} C_{n-3}^{2} C_{n}^{3} - d_{n} \cdot \vert W_{1} \vert + \frac{\vert W_{1} \vert^{2}}{\alpha_{n}}\left(1+h\left(n\right)\right) \geq$$
$$\geq \frac{3}{2} C_{n-3}^{2} C_{n}^{3} - 3 \cdot C_{n-3}^{2} \cdot \vert W_{1} \vert + \frac{\vert W_{1} \vert^{2}}{n} \left(1 + h\left(n\right)\right) =$$
$$ = \frac{3}{2} C_{n-3}^{2} C_{n}^{3} - 3 \cdot C_{n-3}^{2} \cdot c_{n} \cdot C_{n}^{3} + \frac{c_{n}^{2} \cdot \left( C_{n}^{3}\right)^{2}}{n} \left(1 + h\left(n\right)\right) = $$
$$ = \frac{3}{2}C_{n-3}^{2} C_{n}^{3} \left(1 - 2c_{n} + \frac{2}{3}\frac{c_{n}^{2} \cdot C_{n}^{3}}{n \cdot C_{n-3}^{2}} \left(1 + h\left(n\right) \right)\right) = $$
$$ = \frac{n^{5}}{8} \left(1 - \frac{10}{n} + \frac{35}{n^{2}} - \frac{50}{n^{3}} + \frac{24}{n^{4}}\right) \cdot \left(1 - 2c_{n} + \frac{2\cdot c_{n}^{2} \cdot \left( n - 1\right) \left( n - 2\right)}{9 \left(n - 3 \right) \left(n - 4 \right)} \left(1 + h \left( n \right) \right) \right) \geq $$
$$ \geq \frac{n^{5}}{8} \left(1 - \frac{10}{n} \right) \cdot \left(1 - 2c_{n} + \frac{2 \cdot c_{n}^{2}}{9} \left(1 + h\left(n\right) \right) \right) \geq f\left(n\right).$$ 

	
	\item \textit{Функция $l$ удовлетворяет условиям пункта 4 теоремы 3.}
	
	Заметим, что в данном случае множество $W_{1}$ может вовсе не содержать рёбер. Поэтому от улучшения оценки придётся отказаться. 
	Положим 
\begin{equation*}
f(n) = \begin{cases}
0 &\text{\text{при } \; $n \leq 10$}\\
\frac{n^{5}}{8} \left(1 - 2 c_{n} - \frac{10}{n} + \frac{20 \cdot c_{n}}{n}\right) &\text{\text{при } $n > 10$}
\end{cases}
\end{equation*}
Заметим, что $f(n) \sim n^{5} \left(\frac{1}{8} - \frac{c_{n}}{4} \right)$ при $n \rightarrow \infty$. Стало быть, как и в предыдущем пункте, осталось доказать, что $r(l(n)) \geq f(n)$ для любого достаточно большого натурального $n$.
	Для этого запишем следующую цепочку неравенств:
	$$r(l(n)) \geq \vert E(W) \vert \geq \frac{3}{2} C_{n-3}^{2} C_{n}^{3} - d_{n} \cdot \vert W_{1} \vert = \frac{3}{2} C_{n-3}^{2} C_{n}^{3} - 3 \cdot C_{n-3}^{2} \cdot \vert W_{1} \vert =$$
	$$=\frac{n^{5}}{8} \cdot \left(1 - \frac{10}{n} \right) \cdot \left( 1 - 2c_{n}\right) \geq f(n).
	$$ 
\end{enumerate}


\end{document}